\documentclass{amsart}
\usepackage{amsmath}
\usepackage{amsfonts}
\usepackage{amssymb}
\usepackage{graphicx}
\usepackage{float, verbatim}
\usepackage{enumerate}
\usepackage{color}
\usepackage{hyperref}
\usepackage{subcaption}

\newcommand{\Assouad}{\dim_{\mathrm{A}}}

\newtheorem{thm}{Theorem}[section]

\newtheorem{cor}[thm]{Corollary}
\newtheorem{defn}[thm]{Definition}

\newtheorem{rem}[thm]{Remark}

\begin{document}

\title{Assouad dimension of random processes}

\author{Douglas C. Howroyd}
\address{Douglas C. Howroyd\\
	School of Mathematics \& Statistics\\University of St Andrews\\ St Andrews\\ KY16 9SS\\ UK \\ }
\curraddr{}
\email{dch8@st-andrews.ac.uk}
\thanks{}

\author{Han Yu}
\address{Han Yu\\
School of Mathematics \& Statistics\\University of St Andrews\\ St Andrews\\ KY16 9SS\\ UK \\ }
\curraddr{}
\email{hy25@st-andrews.ac.uk}
\thanks{}



\date{}

\dedicatory{Happy Birthday Jon!}

\begin{abstract}
In this paper we study the Assouad dimension of graphs of certain L\'evy processes and functions defined by stochastic integrals. We do this by introducing a convenient condition which guarantees a graph to have full Assouad dimension and then show that graphs of our studied processes satisfy this condition. 
\end{abstract}

\maketitle

\section{Definitions and motivations}\label{LINK}

Studying the dimension of various random processes has been of interest for quite some time. In this paper we will consider the Assouad dimension of graphs of certain random processes, notably L\'evy processes and functions defined by stochastic integrals. 

\emph{L\'{e}vy processes} $X(t)$ were first introduced by Paul L\'evy in 1934 \cite{Le} and are defined to be the stochastic processes satisfying:
\begin{itemize}
	\item[1]: $X(0)=0$ almost surely.
	\item[2]: For all $t,h>0$, $X(t+h)-X(t)$ is equal to $X(h)$ in distribution (stationary increments).
	\item[3]: For all $0<t_1<t_2<...<t_k$ the random variables $X(t_i)-X(t_{i-1})$ are independent (independence of increments).
	\item[4]: For all $t>0$ $\lim_{h\to 0} X(t+h)-X(t)=0$ in probability (continuity).
\end{itemize}
We can construct $X(t)$ such that it is almost surely right continuous with left limits (denoted c\`adl\`ag). Such processes are standard tools in many areas of modern mathematics and its applications. A common example of a L\'evy process is the \emph{Wiener process} (or Brownian motion) where property 2 (stationary increments) is replaced by Gaussian increments, so $X(t+h)-X(t)$ is normally distributed with mean 0 and variance $h$. One can similarly define $d$-dimensional Brownian motion by considering the vector-valued stochastic process $(W_1,\ldots, W_d)$ where the $W_i$ are independent Weiner processes. 

\begin{figure}[h]
\includegraphics[width=0.5\textwidth]{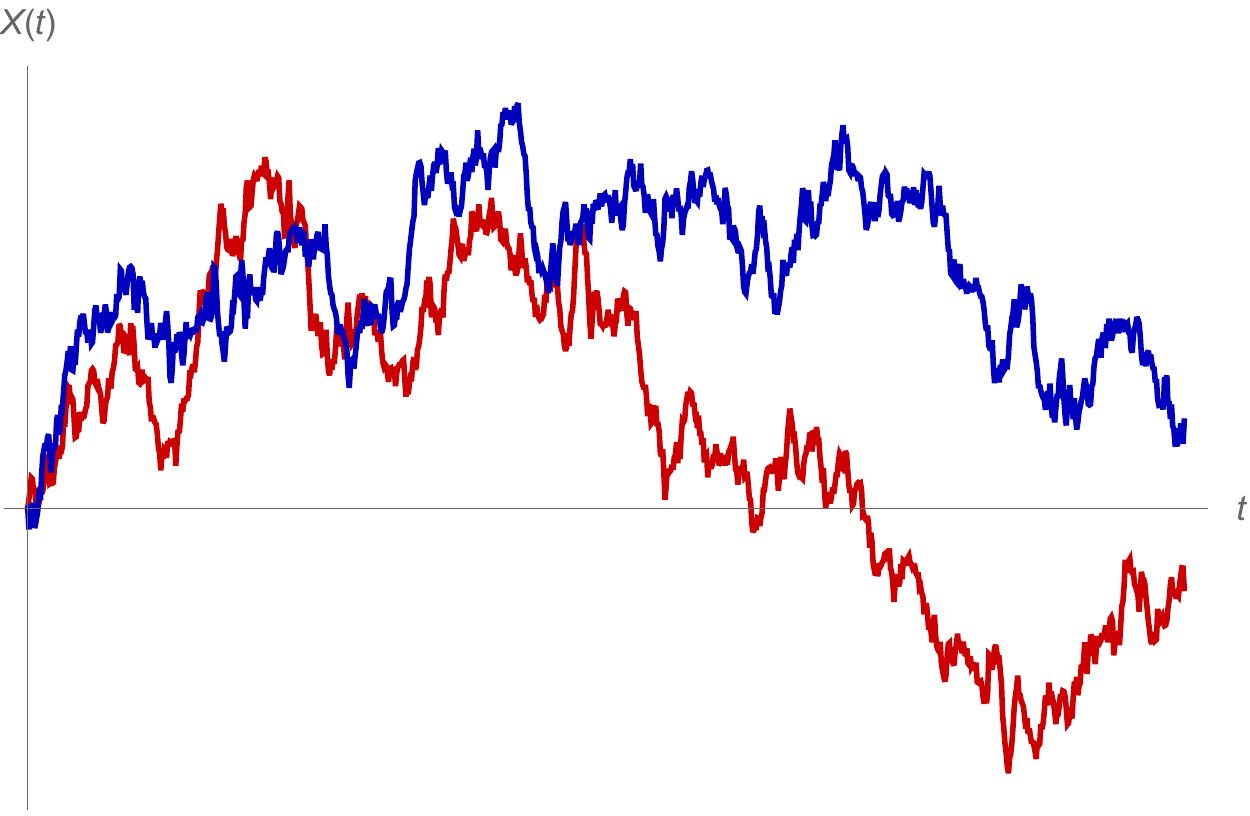}
\caption{\label{fig:brownianmotion}Two graphs of one-dimensional Brownian motion}
\end{figure}

The geometric properties, such as dimension, of Wiener processes have been a particularly well studied area. This includes studying the \emph{graphs}, level sets and \emph{trails} of such processes which can often be thought of as fractals as they often display some \emph{statistical self-affinity}. For any left continuous function $X:\mathbb{R}\to\mathbb{R}$, we define the graph of the function by:
\[
G^{I\subset\mathbb{R}}_{X}=\{(t,y)|y=X(t),t\in I\}\cup J,
\]
where $J$ is the union of vertical segments joining the discontinuities. $J$ is well defined because $X$ is right continuous. It is clear that if $X$ is continuous then $J$ is empty. Taylor \cite{Ta} first calculated the Hausdorff dimension of $d$-dimensional Brownian motion $B_d:\mathbb{R}\to\mathbb{R}^d$ where he showed that almost surely
\[
\dim_\text{H} G_{B_1}^{[0,1]} =  \frac{3}{2}         
\]
and for any $d\ge 2$
\[
\dim_\text{H} B_d([0,1]) =  2.
\]

Another generalisation of Brownian motion is \emph{fractional Brownian motion}, first introduced by Mandelbrot and Van Ness \cite{MVN}. Index-$h$ fractional Brownian motion (fBm) on $\mathbb{R}$ with $0<h<1$ is defined to be the stochastic integral 
\[
B_h (t) = c(h)^{-1} \int_{-\infty}^\infty \left(\left( t-x \right)_+^{h-1/2} - (-x)_+^{h-1/2} \right)W(dx)
\]
where W is the Wiener measure, $(x)_+ = \max\{0,x\}$ and $c(h) = \Gamma(h+1/2)$, where $\Gamma$ denotes the gamma function. Equivalently this is a Gaussian random process $B_h(t)$ with:
\begin{itemize}
	\item[1]: $B_h(0)=0$ almost surely.
	\item[2]: For all $t,u>0$, $B_h(t+u)-B_h(t)$ has normal distribution with mean 0 and variance $u^{2h}$ (Gaussian increments).
	\item[3]: For all $t>0$ $\lim_{u\to 0} B_h(t+u)-B_h(t)=0$ in probability (continuity).
\end{itemize}
One can see that when $h=1/2$, $B_{1/2}(t) = B(t)$ is simply Brownian motion. Much progress has been made on the properties of fBm, see for instance \cite{Ad},\cite{Ka} and \cite{Fa2}. Notably it was shown that almost surely, the graph over the unit interval of index-$h$ fBm has Hausdorff dimension $2-h$. 

\begin{figure}[h]
    \centering
    \begin{subfigure}[b]{0.3\textwidth}
        \includegraphics[width=\textwidth]{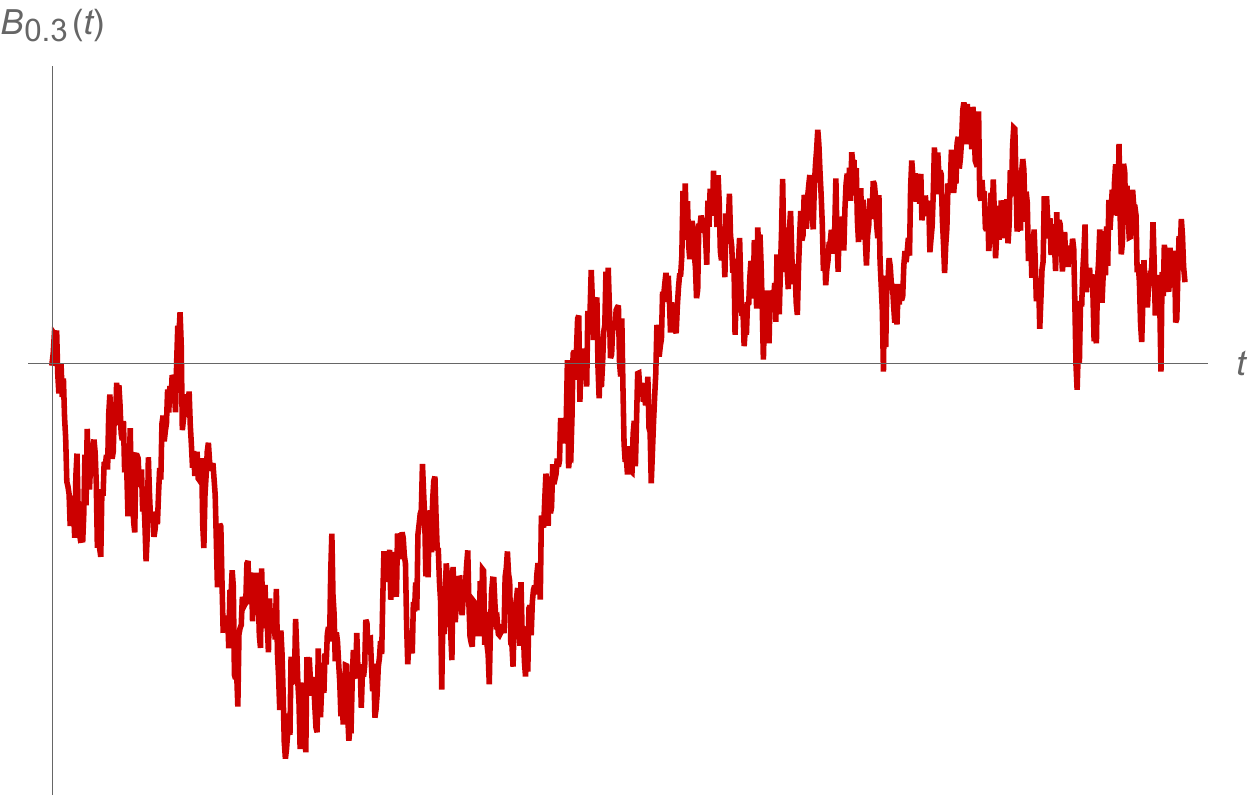}
        \caption{Graph of fBm with h=0.3}
        \label{fig:fbm3}
    \end{subfigure}
    ~ 
    \begin{subfigure}[b]{0.3\textwidth}
        \includegraphics[width=\textwidth]{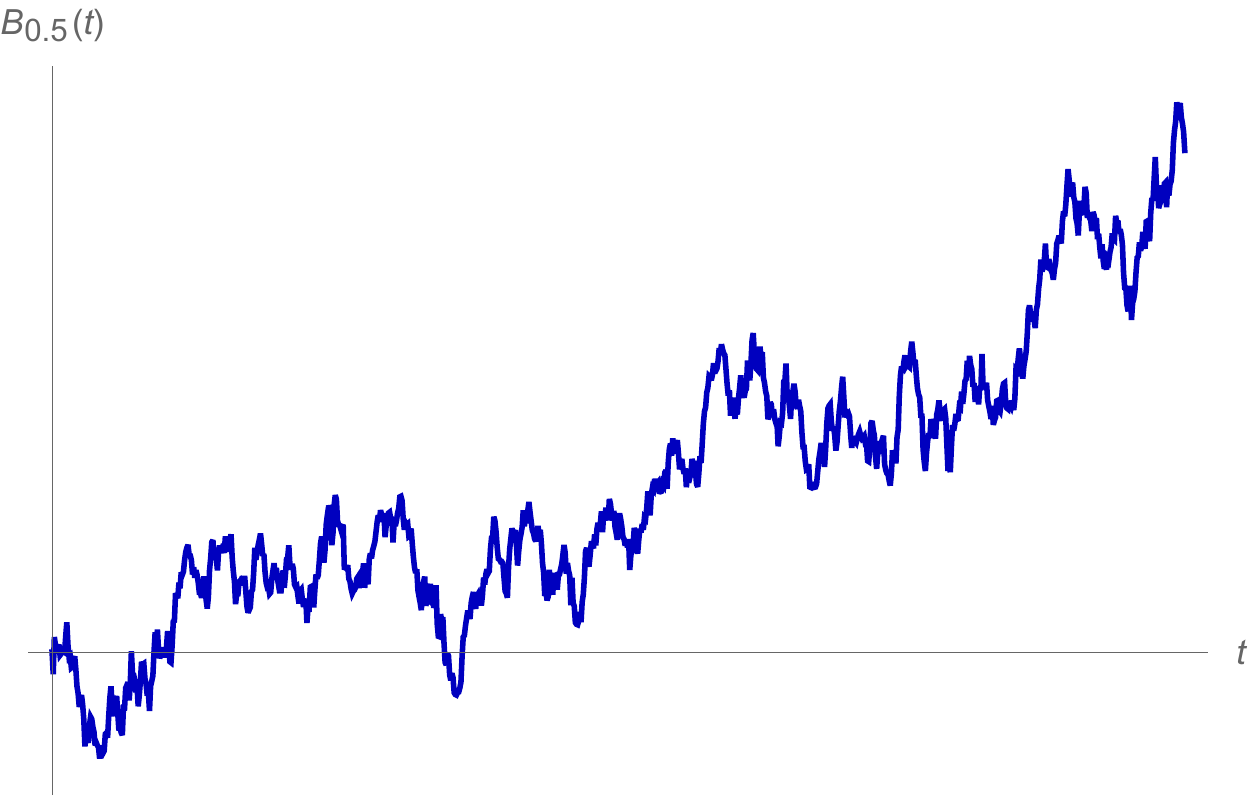}
        \caption{Graph of fBm with h=0.5}
        \label{fig:fbm5}
    \end{subfigure}
    ~ 
    \begin{subfigure}[b]{0.3\textwidth}
        \includegraphics[width=\textwidth]{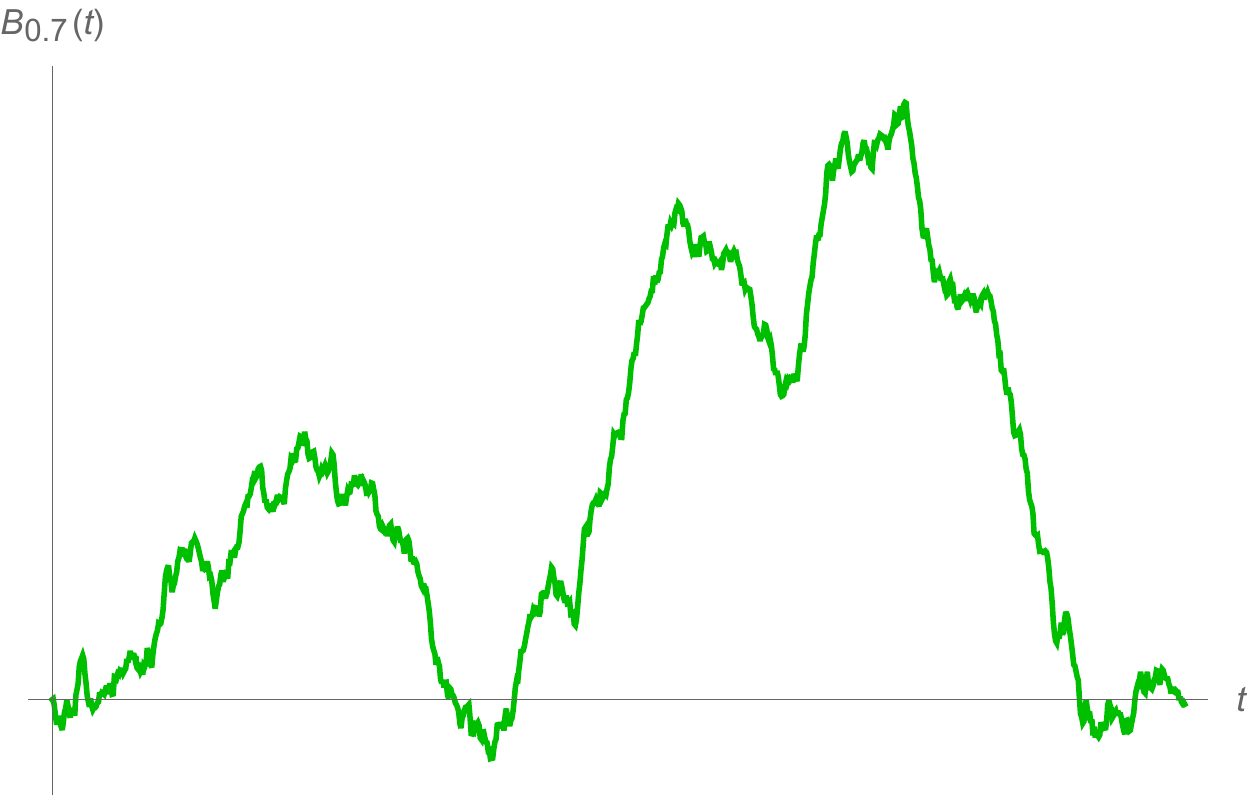}
        \caption{Graph of fBm with h=0.7}
        \label{fig:fbm7}
    \end{subfigure}
\end{figure}

Studying the Assouad dimension of various fractals and its properties is an increasingly popular area of research. In this paper, we are interested in calculating the Assouad dimension of the graph $G_X^{[0,1]}$ for $\beta$-scaling (or L\'{e}vy $\beta$-stable) processes and stochastic integrals $X$. These results will be compared to the previously obtained Hausdorff dimensions.

We say that $X(t)$ satisfies a $\beta$-scaling property if, for any $t,a>0$:
\[
a^{-\frac{1}{\beta}}X(at)=^d X(t),
\]
where `$=^d$' denotes `equal in distribution'.
For example the Wiener process has the $2$-scaling property.
 
A non-empty compact bounded set $F$ is said to be $s$-\textit{homogeneous} if there exists a constant $C>0$ such that for all $0<R$, $r\in (0,R]$ and $x\in F$
\[
N(B(x,R) \cap F, r) \le C\left(\frac{R}{r}\right)^s
\]
where $B(x,R)$ denotes the closed ball of centre $x$ and radius $R$ and $N(E,r)$ is the number of squares of the $r\times r$ grid that intersect the compact set $E$. 

The \emph{Assouad dimension} of a non-empty compact bounded set $F$ is then defined to be 
\[
\Assouad F = \inf \left\{ s \ge 0 \, \, \colon \, F \text{ is } s\text{-homogeneous} \right\}.
\]
This dimension provides information on the extremal local scaling of a set, in this setting it will tell us about the maximal fluctuations of a random process. For a more detailed introduction to this dimension, see \cite{Fr, Ro}. 

This paper will be split into three parts. First we will define a condition, Definition 2.1 which guarantees that a graph will have full Assouad dimension. Then in Section 3 we show that graphs of $\beta$-scaling L\'evy process satisfy this condition and combining these results we prove that graphs of functions defined by certain stochastic integrals also have full Assouad dimension. Finally in Section 4 we remark that our results extend to higher dimensions.

\section{Assouad dimension of graphs}
In this section we will state a convenient condition to check whether a graph of a function $f:[0,1]\to\mathbb{R}$ has full Assouad dimension.

We begin with a definition:

\begin{defn}\label{Win}
	Let $R_1,R_2>0$ be positive numbers and $n_1,n_2>0$ be integers. Given a point $a\in\mathbb{R}^2$, we define $W_{n_1\times n_2}^{R_1\times R_2}(a)$ as the following collection of sets:
	\[	
	\left\{D_{i,j}+a \, \vert \, D_{i,j}=\left[\frac{i}{n_1}R_1,\frac{i+1}{n_1}R_1\right]\times \left[\frac{j}{n_2}R_2,\frac{j+1}{n_2}R_2\right], i\in\{0,\dots,n_1-1\}, j\in\{0,\dots,n_2-1\}\right\}.
	\]
	
	We see that $W_{n_1\times n_2}^{R_1\times R_2}(a)$ is the collection of rectangles with disjoint interiors which partitions the $R_1\times R_2$ rectangle whose bottom left vertex is $a$.
	
	Then let $N_{n_1\times n_2}^{R_1 \times R_2 }(a,G_f^{[0,1]})= \# \left\{ W_{n_1\times n_2}^{R_1\times R_2}(a) \cap G_f^{[0,1]} \right\}$, this is the number of rectangles which intersect the graph.
\end{defn}

The following theorem is a direct consequence of the definition of Assouad dimension and we omit the proof:

\begin{thm}\label{graph}
	If there exists an $A>0$ and sequences:
	\[
	a_i\in \mathbb{R}^2,R_i\in (0,1),\, n_i\in\mathbb{N} \quad (\forall i\in \mathbb{N})
	\]
	with $n_i \rightarrow \infty$ such that for all $i\in\mathbb{N}$
	\[
	N_{n_i\times n_i}^{R_i \times R_i }(a_i,G_f^{[0,1]})\geq A n_i^2.
	\]
	Then
	\[
	\Assouad G_f^{[0,1]}=2.
	\]
\end{thm}

Whilst this might seem like a restrictive condition to ask a general function to satisfy, it is quite natural in the setting of Wiener processes due to the almost sure unbounded variation and $\beta$-scaling property of the process. Considering squares instead of balls in the definition of Assouad dimension is similar to the definition of the Furstenberg star dimension, which is in fact equivalent to the Assouad dimenion, see \cite{chenwuwu}.

\begin{rem}
	Note that one could replace the inequality with the following equality
	\[
	N_{n_i\times n_i}^{R_i \times R_i }(a_i,G_f^{[0,1]})= n_i^2.
	\]
	This follows from \cite[Theorem 2.4]{FY}, where it is shown that a set has full Assouad dimension if and only if it has the unit ball as a weak tangent. This means that any cover of our set is also a cover of a ball and so all smaller squares are needed in the cover.
\end{rem}

\section{Applications to  $\beta$-scaling L\'evy processes}\label{LP}

Let $X(t)$ be a L\'evy process. We assume that $X(1)$ is non-vanishing almost everywhere on $\mathbb{R}$ as a random variable, that is, the distribution function of $X(1)$ is 0 only on a set of measure 0. 

Then for any $\beta$-scaling L\'evy process $X$ we can compute the probability of the following event `$N_{n\times n}^{1 \times 1 }(0,G_X^{[0,1]})=n^2$'. It is a positive number depending only on $n$ and we use $P(n)$ to denote this number. The event `$X$ hits a rectangle' is measurable when $X$ is continuous; when it is discontinuous we join the graph with a vertical line and the process is c\`adl\`ag so the event `$X$ hits a rectangle' is still measurable. Thus our event is measurable as the union of measurable events.

For a $\beta$-scaling random process $X$, we can decompose the graph into countably many disjoint parts:
\[
G_X^{[0,1]}=\bigcup_{i=0}^\infty G_X^{I_i},
\]
where $I_i$ are closed intervals with disjoint interiors such that their union is the unit interval. For our case one could think of this as partitioning the unit interval by intervals of length $1/2^i$. For example take $a_1=0$ and for all $i\geq 1$ let $a_{i+1}=a_i+1/2^i$ and $I_i=[a_i,a_i+1/2^{i}].$

Denote by $|I_i|$ the length of interval $I_i=[a_i,b_i]$. Because we can take $X$ as a left continuous function, $X(a_i)\in\mathbb{R}$ is defined for all $i$. For each $i$ we can apply a linear map $T_i: G_X^{I_i} \rightarrow [0,1]^2$:
\[
T_i(x,y)=\left(\frac{1}{|I_i|}(x-a_i),\frac{1}{|I_i|^{1/\beta}}(y-X(a_i))\right).
\]

By definition, it is easy to see that $T_i(G_X^{I_i})$ are independent $\beta$-scaling L\'evy processes with the same, original distribution. For convenience we identify $T_i(G_X^{I_i})=G_{X_i}^{[0,1]}$ where $X_i$  are independent, identically distributed $\beta$-scaling L\'evy processes.

Let $n_i$ be a sequence of integers such that $\lim_{i\to\infty} n_i=\infty$. We can compute the probability of the event $A_i=`N_{n_i\times n_i}^{1 \times 1 }((0,0),G_{X_i}^{[0,1]})=n_i^2 $'. According to the discussions above, we see that the probability is $P(n_i)$. In fact denote $t_k=k/n^2_i$ for all $k\in\{1,\dots,n^2_i\}$ and $D(j)=[j/n_i,(j+1)/n_i]$ for all $j\in\{0,\dots, n_i-1\}$ then we can see that:
\[
P(n_i)\geq P\left(\forall k\in\left[1,n^2_i\right], k\in \mathbb{N},X(t_k)\in D\left(\left\{\frac{k}{n_i}\right\}n_i\right)\right)>0.
\]
Here $\{ \cdot \}$ denotes the fractional part function. The last inequality follows from our assumption that $X(1)$ is non-vanishing almost everywhere on $\mathbb{R}$. It is clear that this restraint could be relaxed to non-vanishing on some interval without much effort. The rest follows from the property of independent increments.

We can choose $n_i$ to grow slowly enough such that $\sum_{i}P(n_i)=\infty$. Note that the $R_i$ can be chosen so that each square is disjoint and as L\'evy processes are Markov, the events $A_i$ are all independent. Then by Borel-Cantelli lemma we see that with probability $1$, infinitely many events $A_i$ occur. Now if $A_i$ happens, then:
\[
N_{n_i\times n_i}^{1 \times 1 }\left((0,0),G_{X_i}^{[0,1]}\right)=n_i^2,
\]
applying the function $T^{-1}_i$ to the graph, we see that (remember $I_i=[a_i,b_i]$):
\[
N_{n_i\times n_i}^{|I_i| \times |I_i|^{1/\beta} }\left((a_i,X(a_i)),G_{X}^{I_i}\right)=n_i^2.
\]

Since $\beta\geq 1$, we see that $|I_i|\leq |I_i|^{1/\beta}$. Also remember that $X$ can be taken to be a right continuous function and we also include the vertical segments of the jumps in $G^{[0,1]}_X$. Therefore it is clear that there exist an absolute constant $C>0$:
\[
N\left(B\left(\left(a_i,X(a_i)\right),|I_i|\right)\cap G_X^{[0,1]},\frac{|I_i|}{n_i}\right)\geq C n_i^2.
\] 

As infinitely many $A_i$ occur, using Theorem \ref{graph}, we see:
 \[
 \Assouad G_{X}^{[0,1]}=2.
 \]
 
 We conclude the above argument as the following theorem:
 \begin{thm}\label{Main}
 	Let $X$ be a $\beta$-scaling L\'evy process with $\beta\geq 1$, such that $X(1)$ is a random variable whose distribution function is non-vanishing almost everywhere. Then almost surely:
 	\[
 	\Assouad G_{X}^{[0,1]}=2.
 	\]
 \end{thm}

We know that the Wiener process $W(t)$ is a $2$-scaling L\'evy process, therefore we see that:
\begin{cor}
		\[
	\Assouad G_{W}^{[0,1]}=2.
	\]
\end{cor}

Ville Suomala and Changhao Chen, in a personal communication, kindly remarked that this result follows from the graph of Brownian motion having full lower porosity dimension. This approach is inspired by \cite{coxgriffin}, where it was shown that the graph of Brownian motion has full upper porosity dimension. However this porosity dimension technique does not extend to our following, more general result, which relies upon this one.

In fact we can say more about the Assouad dimension of random processes which are functions defined as stochastic integrals, such as fractal Brownian motion. 

\begin{thm}
	Let $f:\mathbb{R}\to\mathbb{R}$ be a function which is zero only finitely often, continuous on some interval and has continuous derivative on that same interval. Then we define $B_f(t)$ as the function defined by the stochastic integral:
	\[	
B_f(t)=\int_0^t f(x) W(dx).
	\]
	We have that almost surely:
	\[
		\Assouad G_{B_f}^{[0,1]}=2.
	\]
\end{thm}
\begin{rem}
	In particular, graphs of fractional Brownian motions with indices $0<h<1$ have full Assouad dimension almost surely.
\end{rem}
\begin{proof}
	Given a function $f$ which is zero only finitely often, continuous and has continuous derivative on some interval, say $J$, we can simply focus on the function restricted to $J$, normalising to obtain a function which is $C^1$ and zero only finitely often on the unit interval. We may then assume that $f(t)>1$ for $t\in [0,1]$ by again restricting our function to an interval where the function is bounded away from zero and normalising. 
    
    As the Assouad dimension provides local information, if the dimension of the graph of the function defined by stochastic integral of this new function is full then the dimension of the original graph is also full. Thus we assume for the rest of this proof that $f$ is a $C^1$ function which is greater than $1$.

Ideally we would wish to integrate by parts in the standard Riemann–-Stieltjes sense
\[
\int_{0}^{t}f(x)W(dx)=f(t)W(t)-\int_{0}^{t}W(x)f'(x)dx.\tag{*}
\]
The problem is that the integral on the left side of the above equation is interpreted as the It\^{o} integral, for which regular integration by parts does not hold. There is however a generalisation of this formula for stochastic integrals which holds as $f$ and $W$ are both semimartingales, see \cite{Pr}[Chapter 2 Section 6] for further details. To be precise we should write the following equation
\[
\int_{0}^{t}f(x)W(dx)=f(t)W(t)-\int_{0}^{t}W(x)f'(x)dx-[f,W]_t.
\]
Here $[f,W]_t$ is the quadratic covariation between $f$ and $W$:

Let $0=t_1<t_2<\dots<t_n=t$ be a partition $P$ of $[0,t]$ and $\vert P \vert$ be the maximum of $t_{k+1}-t_k,k\in\{0,\dots,n-1\}$:
\[
[f,W]_t=\lim_{\vert P\vert\to 0} \sum_{i=1}^{n-1} (f(t_{i+1})-f(t_{i}))(W(t_{i+1})-W(t_i)).
\]
The above convergence is taken in the sense of probability. By using Cauchy-Schwarz we see that:
\[
[f,W]_t\leq [f,f]^{1/2}_t[W,W]^{1/2}_t.
\]
However, it is standard that $[f,f]_t=0$ and $[W,W]_t=t$ as $f$ is $C^1$. So we see that the integral by parts formula (*) is indeed correct for this situation. The integral 
\[
\int_0^{t} W(x)f'(x)dx
\]
is defined to be a random process whose sample space is that of the Wiener process, where fixing a sample path of the Wiener process will determine the integral. We are interested in almost sure properties of this process and will do so by considering almost sure properties of the Wiener process.
  	
	The strategy for the rest of this proof is to choose carefully a typical path of the Wiener process. We denote the sample space of the Wiener process as $\Omega$.

	First, we see that for almost all $\omega\in\Omega$, $W(t,\omega)$ is c\`adl\`ag in $t$, and therefore there is a constant $C_\omega$ such that:
	\[
	|W(t,\omega)|\leq C_\omega
	\]
    for all $t\in [0,1]$. 

The second almost sure property is described in the proof of Theorem \ref{Main}, that there are infinitely many intervals $I_i=[a_i,b_i]\subset [0,1]$ and a sequence $n_i\to\infty$ such that for $k\in\{0,1,\dots,n_i-1\}$
    \[
    \left|W\left(a_i+(k+1)\frac{|I_i|}{n_i},\omega\right)-W\left(a_i+k\frac{|I_i|}{n_i},\omega\right)\right|\geq |I_i|^{1/2}\geq |I_i|.
    \]

    In the following discussion we shall fix a typical $\omega$ such that $W(t,\omega)$ satisfies the above two almost sure properties, in particular, we think of $C_\omega>0$ as a fixed constant.

	Then we see that:
	\begin{eqnarray*}
	& &\left|B_f\left(a_i+(k+1)\frac{|I_i|}{n_i},\omega\right)-B_f\left(a_i+k\frac{|I_i|}{n_i},\omega\right)\right|\\
    &=& \left| \int_{a_i+k\frac{|I_i|}{n_i}}^{a_i+(k+1)\frac{|I_i|}{n_i}} f(x)W(dx)\right|\\
	&=& \left| f\left(a_i+k\frac{|I_i|}{n_i}\right)W(x,\omega)\bigg\vert_{a_i+k\frac{|I_i|}{n_i}}^{a_i+(k+1)\frac{|I_i|}{n_i}} + \int_{a_i+k\frac{|I_i|}{n_i}}^{a_i+(k+1)\frac{|I_i|}{n_i}} W(x)f'(x)dx \right|.
	\end{eqnarray*}
	 
	Since $f$ is $C^1$, we see that there is a constant $C_f$ (which does not depend on $i$) such that for all $x\in \left[{a_i+k\frac{|I_i|}{n_i}},{a_i+(k+1)\frac{|I_i|}{n_i}}\right]$:
	\[
	\left\vert W(x)f'(x)\right\vert \leq C_f C_\omega.
	\]
	
	Then we have the following inequalities:
	\[
	\left| f\left(a_i+k\frac{|I_i|}{n_i}\right)W(x,\omega)\bigg\vert_{a_i+k\frac{|I_i|}{n_i}}^{a_i+(k+1)\frac{|I_i|}{n_i}} \right| \geq \left| f\left(a_i+k\frac{|I_i|}{n_i}\right) \right| |I_i|^{1/2}\geq |I_i|^{1/2}.\tag{**}
	\]
	and
	\[
	\left\vert\int_{a_i+k\frac{|I_i|}{n_i}}^{a_i+(k+1)\frac{|I_i|}{n_i}} W(x)f'(x)dx \right\vert\leq C_fC_\omega\frac{|I_i|}{n_i}.
	\]
       Since $|I_i|\to 0$,  we see that for a constant $C'$ which depends only on $f$ and $\omega$:
	\[
	\left|B_f\left(a_i+(k+1)\frac{|I_i|}{n_i},\omega\right)-B_f\left(a_i+k\frac{|I_i|}{n_i},\omega\right)\right|\geq C'|I_i|.\tag{***}
	\]
	 The above inequality holds for all $k\in\{0,\dots,n_i-1\}$. 
     
      Moreover, $W$ has a `zigzag' property. For even integers $k$ we have
\[
W\left(a_i+(k+1)\frac{|I_i|}{n_i},\omega\right)-W\left(a_i+k\frac{|I_i|}{n_i},\omega\right)>0,
\] 
for odd integers $k$ we have
\[
W\left(a_i+(k+1)\frac{|I_i|}{n_i},\omega\right)-W\left(a_i+k\frac{|I_i|}{n_i},\omega\right)<0.
\]  
Heuristically this says that the process increases on the first interval, decreases on the second and so forth, zigzagging from top to bottom. We can see that the expressions inside the absolute values in (**) and (***) also satisfy similar `zigzag' properties. Therefore there is a constant $A=A(\omega,f)>0$ such that:
\[
N\left(B\left(\left(a_i,B_f(a_i)\right),|I_i|\right)\cap G_{B_f}^{[0,1]},\frac{|I_i|}{n_i}\right)\geq A n_i^2.
\] 

This concludes the proof because the above argument holds for a set of full probability $\omega\in\Omega$.
\end{proof}

\section{A remark on higher dimensional Brownian motion}

Definition \ref{Win} has a natural generalization in $\mathbb{R}^d$.
		\begin{defn}\label{Win2}
		Let $R_1,\dots,R_d>0$ be positive numbers and $n_1,\dots, n_d$ be integers. Given a point $a\in\mathbb{R}^d$, we define $W_{n_1\times\dots\times n_d}^{R_1\times\dots\times R_d}(a)$ as the following collection of sets:
		\begin{eqnarray*}
		\Bigg\{D_{i_1,\dots,i_d}+a \, \vert \, D_{i_1,\dots,i_d}=\left[\frac{i_1}{n_1}R_1,\frac{i_1+1}{n_1}R_1\right]\times\dots\times \left[\frac{i_d}{n_d}R_d,\frac{i_d+1}{n_d}R_d\right],\\ i_j\in\{0,\dots,n_j-1\}, j\in\{1,\dots,d\}\Bigg\}.
		\end{eqnarray*}
		
		We see that $W_{n_1\times n_2}^{R_1\times R_2}(a)$ is a collection of rectangles with disjoint interiors.
	\end{defn}

Using the above definition and a similar argument as the one in section \ref{LP},  Theorem \ref{graph} also extends to higher dimensions, we can show the following result and we omit the proof:

\begin{thm}
	Let $B_d(t)$ be the $d$ dimensional Brownian motion, from $\mathbb{R}$ to $\mathbb{R}^d$. Then almost surely:
	\[
	\Assouad B_d([0,1])=d.
	\]
\end{thm}

We can compare this result to the well known \[\dim_{\mathrm{H}} B_d([0,1])=2\] for $d\geq 2$. The Hausdorff dimension being 2 here can be thought as a reflection that higher dimensional Brownian motion is transient, whilst the Assouad dimension shows that there are still areas of maximal fluctuation. 

 Brownian motion also provides examples of Salem sets that can have different Hausdorff and Assouad dimensions, we refer the reader to \cite{Ka}  for further discussion on the links between random processes and Salem sets.

\begin{centering}

\textbf{Acknowledgements}

\end{centering}

DCH was financially supported by EPSRC Doctoral Training Grant EP/N509759/1. HY thanks the University of St Andrews' School of Mathematics and Statistics for the PhD scholarship.

\providecommand{\bysame}{\leavevmode\hbox to3em{\hrulefill}\thinspace}
\providecommand{\MR}{\relax\ifhmode\unskip\space\fi MR }
\providecommand{\MRhref}[2]{%
	\href{http://www.ams.org/mathscinet-getitem?mr=#1}{#2}
}
\providecommand{\href}[2]{#2}

\end{document}